\documentclass[11pt,reqno]{amsart}
\usepackage{amsmath,amssymb,amsthm}

 \newcommand{\dd}{{\mathrm{d}}} 
 \newcommand{\ii}{{\mathrm{i}}} 
 \newcommand{\I}{1\!\!1}    

\newcommand{\Prob}{{\mathbb{P}}}
 \DeclareMathOperator{\E}{{\mathbb{E}}}

\newcommand{\C}{{\mathbb{C}}}

\renewcommand{\Im}{{\mathfrak{Im}}}

\renewcommand{\t}{\theta}

\newtheorem{thm}{Theorem}
 
 \newtheorem{lem}[thm]{Lemma}

 \theoremstyle{remark}
 \newtheorem*{rem}{Remark}

\begin{document}
\title{An arithmetic model for the total disorder process}
\author{C.P. Hughes}
\address{Department of Mathematics, University of York, Heslington,
York, YO10 5DD, U.K.}
 \email{ch540@york.ac.uk}

\author{A. Nikeghbali}
 \address{ETHZ \\ Departement Mathematik, R\"{a}mistrasse 101, HG G16\\ Z\"{u}rich 8092, Switzerland}
 \email{ashkan.nikeghbali@math.ethz.ch}

\author{M. Yor}
 \address{Laboratoire de Probabilit\'es et Mod\'eles Al\'eatoires \\
Universit\'e Pierre et Marie Curie, et C.N.R.S. UMR 7599 \\ 175,
rue du Chevaleret \\ F-75013 Paris, France}

\date{1 December 2006}

\begin{abstract}
We prove a multidimensional extension of Selberg's central limit
theorem for the logarithm of the Riemann zeta function on the
critical line. The limit is a totally disordered process, whose
coordinates are all independent and Gaussian.
\end{abstract}

\maketitle

\section{Introduction}

A classical result of Selberg \cite{Selb89} (see also Laurin{\v
c}ikas, \cite{Laurin}) states that the classical continuous
determination of the logarithm of the Riemann zeta function is
asymptotically normally distributed, in the sense that if $\Gamma$
is a regular Borel measurable subset of $\C$,
\begin{equation*}
\lim_{T\to\infty} \frac{1}{T} \int_T^{2T} \I\left\{
\frac{\log\zeta(\tfrac12 + \ii t)}{\sqrt{\tfrac12 \log\log T}}\in
\Gamma \right\} \;\dd t = \frac{1}{2\pi} \int_{\Gamma}
e^{-(x^2+y^2)/2} \;\dd x \;\dd y
\end{equation*}
where $\I$ is the indicator function, and regular means that the
boundary of $\Gamma$ has zero Lebesgue measure.

If we let
\begin{equation*}
L_\lambda(N,u) := \frac{\log\zeta(\tfrac12 + \ii u
e^{N^\lambda})}{\sqrt{\log N}}
\end{equation*}
then Selberg's result implies that
\begin{equation*}
\lim_{N\to\infty} \int_1^2 \I\left\{ L_\lambda(N,u) \in \Gamma
\right\} \;\dd u = \Prob\{G_{\lambda}\in\Gamma\}
\end{equation*}
where $G_\lambda = G_\lambda^{(1)} + \ii G_\lambda^{(2)}$ is a
complex-valued Gaussian random variable with mean zero and
variance $\lambda/2$, i.e.: $G_\lambda^{(1)}$ and
$G_\lambda^{(2)}$ are independent, centered, and
$\E[(G_\lambda^{(1)})^2] = \E[(G_\lambda^{(2)})^2] = \lambda/2$.

It is now a natural question, at least from a probabilistic
standpoint, to look for an asymptotic distribution for
$\left(L_{\lambda_1}(N,\cdot) , \dots
,L_{\lambda_k}(N,\cdot)\right)$, for different $\lambda_{i}$'s.

\begin{thm}\label{thm:main}
For $\lambda_1>\lambda_2>\dots>\lambda_k>0$, and for every
$\left(\Gamma_{i},i\leq k\right)$ regular,
\begin{equation}\label{eq:main}
\lim_{N\to\infty} \int_1^2 \I\left\{L_{\lambda_1}(N,u) \in
\Gamma_1 , \dots , L_{\lambda_k}(N,u) \in \Gamma_k \right\} \;\dd
u = \prod_{j=1}^k \Prob\left\{G_{\lambda_j} \in \Gamma_j \right\}
.
\end{equation}
\end{thm}

We now note that if $\left(D_\lambda = D_\lambda^{(1)} + \ii
D_\lambda^{(2)} , \lambda>0\right)$ is a totally disordered
complex-valued Gaussian process, meaning that
$\left(D_\lambda^{(1)}
, \lambda>0\right)$ and \\
 $\left(D_\lambda^{(2)} ,
\lambda>0\right)$ are two independent Gaussian processes all of
whose coordinates are independent with $\E[(D_\lambda^{(1)})^2] =
\E[(D_\lambda^{(2)})^2] = \lambda/2$, then the quantity on the right
hand side of \eqref{eq:main} is
\[
\Prob\left\{D_{\lambda_1} \in \Gamma_1 , \dots , D_{\lambda_k} \in
\Gamma_k\right\}.
\]

Theorem \ref{thm:main} is an attempt to move from the
deterministic set up of the Riemann zeta function, and the
``static'' central limit theorem of Selberg into a more
``dynamic'' probabilistic world, where a process appears in the
limit. However, this process is quite wild. In the next section,
we comment about it, and some of its occurrences in random matrix
theory. Finally, in the third section we prove
Theorem~\ref{thm:main} using the method of moments.

\begin{rem}
Our methods apply equally well to any $L$-function from the Selberg
class, but for concreteness and for the sake of simplicity we only
state here the result for the Riemann zeta function.
\end{rem}

\section{Some remarks on total disorder process}

\subsection{Non-measurability of the total disorder process}
The total disorder process is a ``wild'' process; indeed there is no
measurable process $(\lambda,\omega) \mapsto \tilde
D_\lambda(\omega)$ which would be a modification of $(D_\lambda ,
\lambda \geq 0)$, i.e. $\Prob\{ \tilde D_\lambda = D_\lambda) = 1$
for all $\lambda$. Indeed, if so, we would get (use Fubini)
$$\int_{a}^{b}\widetilde{D}_{\lambda}\dd \lambda=0\quad \mathrm{a.s.},$$hence$$\widetilde{D}_{\lambda}=0\quad \dd\lambda \dd\mathbb{P},$$which is absurd (for some further discussion on the
total disorder process, see page 37 of \cite{RevYor}).

\subsection{The total disorder process in random matrix theory}

The total disorder process has already been observed
asymptotically in random matrix theory, although in a different
guise. Let $Z_U(\t) = \det (I-U e^{-\ii\t})$ be the characteristic
polynomial of an $N\times N$ unitary matrix $U$ chosen with Haar
measure, then Hughes, Keating and O'Connell \cite{HKO2} prove that
\begin{equation*}
\frac{\log Z_U(\theta)}{\sqrt{\tfrac12\log N}}
\end{equation*}
weakly converges to $X(\t)+\ii Y(\t)$, where $X(\t),Y(\t)$ are
independent Gaussian processes   with covariance structure
\begin{equation*}
\E \left[ X(\t_1) X(\t_2) \right] = \E \left[ Y(\t_1) Y(\t_2)
\right] =
\begin{cases}
1 & \text{ if } \t_1 = \t_2\\
0 & \text{ otherwise}
\end{cases}
\end{equation*}

This was used to provide an explanation for the covariance structure
of $C_U(s,t)$, the number of eigenangles of $U$ that lie in the
interval $(s,t)$, found earlier by Wieand \cite{Wieand,Wieand1}. A
separate explanation was given by Diaconis and Evans \cite{DE}. Let
\begin{equation*}
\widetilde C_U(s,t) := \frac{C_U(s,t)-
(t-s)N/2\pi}{\frac{1}{\pi}\sqrt{\log N}} .
\end{equation*}
Wieand  proves that for fixed $s,t$, if the matrices $U$ are chosen
with Haar measure from the unitary group, then $\widetilde C_U(s,t)$
converges in distribution, as $N\to\infty$, to a standard normal
random variable. In fact she goes much further by proving weak
convergence of $\widetilde C_U(s,t)$ to a certain Gaussian process
$C(s,t)$.

\begin{thm}[Wieand]\label{thm:wieand} For
$-\pi<s<t\leq \pi$, the finite dimensional distributions of the
process $\widetilde C_U(s,t)$ converge as $N\to\infty$ to those of
a centered Gaussian process $C(s,t)$ with covariance structure

\begin{equation*}
\E\left\{C(s,t) C(s',t')\right\} =
\begin{cases}
1 & \text{ if } s=s', t=t'\\
-1 & \text{ if } s=t', t=s'\\
\tfrac{1}{2} & \text{ if } s=s' \text{ or if } t=t' \text{ but not both}\\
-\tfrac{1}{2} & \text{ if } s=t' \text{ or if } t=s' \text{ but not both}\\
0 & \text{ otherwise}
\end{cases}
\end{equation*}
\end{thm}
A similar process result had previously been found by Costin and
Lebowitz \cite{CL} for GUE matrices, and Soshnikov \cite{Sosh}
considers a process result for counting the number of eigenangles in
an interval with a given minimum displacement. The surprising thing
about these correlations is that they imply that even if an interval
$I$ contains more than the average number of eigenangles, then any
subset of $I$ not sharing a common endpoint with $I$ will usually
still contain \emph{its} average number. Also, no matter how close
 two intervals $I$ and $J$ are, unless they share an endpoint, then $C_{I}$ and $C_{J}$ (with obvious notations) are independent.

Of course, the results of Wieand and of Hughes, Keating and
O'Connell are strongly related, because
\begin{equation*}
\widetilde C_U(s,t) =  \frac{\Im\log Z_U(t)}{\sqrt{\log N}} -
\frac{\Im\log Z_U(s)}{\sqrt{\log N}} .
\end{equation*}

\section{Proof of theorem~\ref{thm:main}}

We will prove Theorem~\ref{thm:main} via the method of moments.
The following lemma will play an essential role in our argument.

\begin{lem}\label{lem:little_lemma}
A complex random variable $Z=X+\ii Y$ has moments
\begin{equation}\label{eq:little_lemma}
\E\left[ Z^m \overline{Z}^n \right] =
\begin{cases}
n! 2^n \sigma^{2n} & \text{ if } m=n\\
0 & \text{ otherwise}
\end{cases}
\end{equation}
if and only if $X$ and $Y$ are independent and distributed
according to the normal law with mean zero and variance
$\sigma^2$.
\end{lem}

\begin{proof}
Let $Z=X+\ii Y$ where $X$ and $Y$ are independent, centered normal
random variables with variance $\sigma^2$. Consider the joint
moment generating function of $Z$ and $\overline Z$: for
$\left(\alpha,
\beta\right)\in\mathbb{C}^{2}$,
\begin{align*}
\E\left[ e^{\alpha Z} e^{\beta \overline Z}\right] &=
\E\left[e^{(\alpha+\beta)X + \ii(\alpha-\beta)Y}\right] \\
&=e^{(\alpha+\beta)^2\sigma^2/2 - (\alpha-\beta)^2\sigma^2/2} \\
&= e^{2\alpha\beta\sigma^2} \\
&= \sum_{n=0}^\infty n! 2^n \sigma^{2n} \frac{\alpha^n}{n!}
\frac{\beta^n}{n!}
\end{align*}
which is the two-variable moment generating function of
\eqref{eq:little_lemma}. Conversely, assume that
\eqref{eq:little_lemma} holds for the joint moments of $Z$ and
$\overline Z$. Then working up the above chain of equalities
proves that $Z=X+\ii Y$ where $X$ and $Y$ are independent,
centered gaussians with variance $\sigma^2$.
\end{proof}

From Lemma~\ref{lem:little_lemma}, if one can show that for any
positive integers $k$ and any integers $m_1,\dots,m_k ;
n_1,\dots,n_k$, and if for any $\lambda_1>\dots>\lambda_k$
\begin{equation*}
\E\left[ \prod_{\ell=1}^k D_{\lambda_\ell}^{m_\ell}
\overline{D_{\lambda_\ell}}^{n_\ell} \right] = \prod_{\ell=1}^k
n_\ell! \lambda_\ell^{n_\ell} \delta_{m_\ell , n_\ell}
\end{equation*}
then one may conclude that $\left(D_\lambda , \lambda>0\right)$ is
a centered complex-valued Gaussian totally disordered process with
covariance structure
\begin{equation*}
\E\left[ D_{\lambda_i} \overline{D_{\lambda_j}} \right] =
\begin{cases}
\lambda_i & \text{ if } \lambda_i=\lambda_j\\
0 & \text{ otherwise}
\end{cases}
\end{equation*}
and $\E\left[ D_{\lambda_i} D_{\lambda_j} \right] = 0 $ for all
$\lambda_i, \lambda_j$.

Therefore, Theorem~\ref{thm:main} is a consequence of the
following
\begin{thm}\label{lem:mmts_zeta}
Let
\begin{equation*}
L_\lambda(N,u) = \frac{\log\zeta(\tfrac12 + \ii u
e^{N^\lambda})}{\sqrt{\log N}} .
\end{equation*}
If $\lambda_1<\lambda_2<\dots<\lambda_k$ are fixed, then
\begin{equation*}
\lim_{N\to\infty} \int_1^2 \prod_{j=1}^k L_{\lambda_j}(N,u)^{m_j}
\overline{L_{\lambda_j}(N,u)}^{n_j} \;\dd u = \prod_{j=1}^k n_j!
\lambda_j^{n_j} \delta(m_j , n_j)
\end{equation*}
where $\delta(m_j,n_j) = 1$ if $m_j=n_j$ and zero otherwise.
\end{thm}

We need the following theorem of Selberg, \cite{Selb89}.
\begin{thm}[Selberg]\label{thm:Selberg}
If $n$ is a positive integer, $0<a<1$, and $T^{a/n} \leq x \leq
T^{1/n}$, then
\begin{equation*}
\frac{1}{T} \int_T^{2T} \left| \log\zeta(\tfrac12+\ii t) -
\sum_{p\leq x} \frac{p^{-\ii t}}{\sqrt{p}} \right|^{2n} \;\dd t =
O\left(n^{4n} e^{An} \right)
\end{equation*}
for some constant $A$ which depends upon $a$.
\end{thm}

We also need to calculate the moments of certain prime sums.
\begin{lem}\label{lem:mmts_primes}
Given $k$ a positive integer, let $\lambda_1 > \dots >
\lambda_k>0$. Let
\begin{equation*}
P(\lambda,n) = P(\lambda,n;k,N,u) = \frac{1}{\sqrt{\log N}}
\sum_{p \leq \exp\left(\frac{N^{\lambda}}{40 k n}\right)}
\frac{p^{-\ii u e^{N^\lambda}}}{\sqrt{p}}
\end{equation*}
For any non-negative integers $m_1,\dots,m_k$ and $n_1,\dots,n_k$,
\begin{equation}\label{eq:mmts_primes_show_this}
\lim_{N\to\infty} \int_1^2 \prod_{j=1}^k P(\lambda_j,m_j)^{m_j}
\overline{P(\lambda_j,n_j) }^{n_j} \;\dd u  = \prod_{j=1}^k (n_j)!
\left(\lambda_j\right)^{n_j} \delta(m_j,n_j)
\end{equation}
\end{lem}
(Note that for the sake of simplicity the variable $u$ does not
appear explicitly in (\ref{eq:mmts_primes_show_this}) and in some
expressions below)

\begin{proof}
We wish to expand out
\begin{equation*}
\prod_{j=1}^k P(\lambda_j,m_j)^{m_j}
\end{equation*}
as a multiple sum over primes. It is exceedingly complicated. We
will introduce the following notation: For $j=1,\dots,k$, let
\begin{equation*}
\mathbf{p}_j = \left(p_{j,1},\dots,p_{j,m_j}\right)
\end{equation*}
and let
\begin{equation*}
\mathcal{P}_j =\mathcal{P}\left(\lambda_{j}, m_{j}\right)= \left\{
\mathbf{p}_j \ : \ p_{j,\ell} \text{ is prime } \ , \ p_{j,\ell}
\leq \exp\left(\frac{N^{\lambda_j}}{40 k m_j}\right) \right\}
\end{equation*}
Hence
\begin{multline*}
\prod_{j=1}^k P(\lambda_j,m_j)^{m_j} =(\log N)^{-(m_1 + \dots + m_k)/2} \\
\times \sum_{\mathbf{p}_1 \in \mathcal{P}_1 , \dots , \mathbf{p}_k
\in \mathcal{P}_k} \frac{\exp\left(-\ii u \sum_{j=1}^k
e^{N^{\lambda_j}} \log(p_{j,1}\dots p_{j,m_j})
\right)}{\sqrt{\prod_{j=1}^k \prod_{\ell_j=1}^{m_j} p_{j,\ell_j}}}
\end{multline*}
Similarly, we let
\begin{equation*}
\mathbf{q}_j = \left( q_{j,1} , \dots , q_{j,n_j} \right)
\end{equation*}
and let
\begin{equation*}
\mathcal{Q}_j \equiv\mathcal{P}\left(\lambda_{j},n_{j}\right)=
\left\{ \mathbf{q}_j \ : \ q_{j,\ell} \text{ is prime } \ , \
q_{j,\ell} \leq \exp\left(\frac{N^{\lambda_j}}{40 k n_j}\right)
\right\}
\end{equation*}
and so
\begin{multline*}
\prod_{j=1}^k \overline{P(\lambda_j,n_j)}^{n_j} = (\log N)^{-(n_1 + \dots + n_k)/2} \\
\times \sum_{\mathbf{q}_1 \in \mathcal{Q}_1 , \dots , \mathbf{q}_k
\in \mathcal{Q}_k} \frac{\exp\left(\ii u \sum_{j=1}^k
e^{N^{\lambda_j}} \log(q_{j,1}\dots q_{j,n_j})
\right)}{\sqrt{\prod_{j=1}^k \prod_{\ell_j=1}^{n_j} q_{j,\ell_j}}}
\end{multline*}
Finally, let $\mathbf{p} = \bigcup_{j=1}^k \mathbf{p}_j$ and
$\mathbf{q} = \bigcup_{j=1}^k \mathbf{q}_j$, and let
\begin{equation*}
F(\mathbf{p} , \mathbf{q}) = \sum_{j=1}^k
\exp\left(N^{\lambda_j}\right) \log\left(\frac{q_{j,1}\dots
q_{j,n_j}}{p_{j,1}\dots p_{j,m_j}}\right)
\end{equation*}
Therefore,
\begin{multline}\label{eq:prod_expanded}
\prod_{j=1}^k P(\lambda_j,m_j)^{m_j}
\overline{P(\lambda_j,n_j) }^{n_j}  = (\log N)^{-(m_1 + \dots + m_k + n_1 + \dots + n_k)/2}\\
\times  \sum_{\substack{\mathbf{p}_1 \in \mathcal{P}_1 , \dots ,
\mathbf{p}_k \in \mathcal{P}_k \\ \mathbf{q}_1 \in \mathcal{Q}_1 ,
\dots , \mathbf{q}_k \in \mathcal{Q}_k}} \frac{\exp\left(\ii u
F(\mathbf{p} , \mathbf{q}) \right)}{\sqrt{\prod_{j=1}^k \left(
\prod_{\ell_j=1}^{m_j} p_{j,\ell_j} \right)
\left(\prod_{\ell_j=1}^{n_j} q_{j,\ell_j}\right)}}
\end{multline}
We divide the sum up into two parts, depending on whether
$F(\mathbf{p} , \mathbf{q})$ equals zero or not. The terms where
the sum vanishes we call \emph{diagonal terms}; the other terms
are \emph{off-diagonal}. The proof of the lemma will follow from
showing that the off-diagonal terms do not contribute in the
large-$N$ limit, and using a simple combinatorial enumeration of
the diagonal terms, along with the prime number theorem, to
estimate the diagonal terms.

\subsection{The diagonal terms}\label{sect:diag}

We will see below in \S\ref{sect:off-diag} that since
$\lambda_1>\dots>\lambda_k$, then for sufficiently large $N$, the
only way for $F(\mathbf{p},\mathbf{q})=0$ is if
\begin{equation*}
\exp\left(N^{\lambda_j}\right) \log\left(\frac{q_{j,1}\dots
q_{j,n_j}}{p_{j,1}\dots p_{j,m_j}}\right) = 0
\end{equation*}
for each $j=1,\dots,k$ separately. Thus the diagonal terms are
those contained in the sets
\begin{equation}\label{eq:defn_D}
\mathcal{D}_j := \left\{ (\mathbf{p}_j , \mathbf{q}_j) \ : \
\mathbf{p}_j \in \mathcal{P}_j , \mathbf{q}_j \in \mathcal{Q}_j ,
\prod_{\ell=1}^{m_j} p_{j,\ell} = \prod_{\ell=1}^{n_j} q_{j,\ell}
\right\}
\end{equation}
Since $p_{j,\ell}$ and $q_{j,\ell}$ are both prime, the set
$\mathcal{D}_j$ is empty unless $m_j = n_j$. Under such an
assumption, the diagonal terms in \eqref{eq:prod_expanded} are
\begin{equation*}
\sum_{(\mathbf{p}_1,\mathbf{q}_1) \in \mathcal{D}_1 , \dots ,
(\mathbf{p}_k , \mathbf{q}_k) \in \mathcal{D}_k}
\frac{1}{\prod_{j=1}^k \prod_{\ell_j=1}^{n_j} q_{j,\ell_j}} =
\prod_{j=1}^k \left(\sum_{(\mathbf{p}_j,\mathbf{q}_j) \in
\mathcal{D}_j} \frac{1}{q_{j,1} \dots q_{j,n_j}} \right)
\end{equation*}
If $q_{j,1}, \dots , q_{j,n_j}$ are distinct primes, then there
are $(n_j)!$ ways of choosing $p_{j,1} , \dots, p_{j,n_j}$ such
that the products are equal. (If the $q_{j,\ell}$ are not
distinct, then the result is similar, but with a different
combinatorial factor, and the result is at least a couple of
logarithms smaller). Hence
\begin{align}\label{eq:diag_terms}
\sum_{(\mathbf{p}_j,\mathbf{q}_j) \in \mathcal{D}_j}
\frac{1}{q_{j,1} \dots q_{j,n_j}} &= (n_j)! \sum_{\mathbf{q}_j \in
\mathcal{Q}_j} \frac{1}{q_{j,1} \dots q_{j,n_j}} \left(1+O(\frac{1}{\log^2 N})\right)\nonumber \\
& = (n_j)! \left(\sum_{q\leq \exp\left(\frac{N^{\lambda_j}}{40 k
n_j}\right)}
\frac{1}{q}\right)^{n_j} \left(1+O(\frac{1}{\log^2 N})\right) \nonumber \\
&= (n_j)! \left(\log
\left(\frac{N^{\lambda_j}}{40 k n_j}\right)+O(1)\right)^{n_j} \left(1+O(\frac{1}{\log^2 N})\right) \nonumber \\
& = (n_j)! (\lambda_j \log N)^{n_j} \left(1+O(\frac{1}{\log
N})\right)
\end{align}

Hence the diagonal contribution to \eqref{eq:prod_expanded} is
\begin{equation*}
\left(\prod_{j=1}^k (n_j)! (\lambda_j)^{n_j} \delta(m_j ,
n_j)\right)\left(1+O(\frac{1}{\log N})\right)
\end{equation*}
which is the right-hand side of \eqref{eq:mmts_primes_show_this}
in the large-$N$ limit. (The constant implicit in the $O$-term
depends on $m_j, n_j, \lambda_j$ and $k$, but these are all
constants). Hence the proof of the lemma will be complete if we
can show there is no contribution to
\eqref{eq:mmts_primes_show_this} from the off-diagonal terms, the
terms where $F(\mathbf{p} , \mathbf{q}) \neq 0$.

\subsection{The off-diagonal terms}\label{sect:off-diag}

Now we  show that the non-diagonal terms of
\eqref{eq:prod_expanded} do not contribute to
\eqref{eq:mmts_primes_show_this} in the limit. Upon integrating
\eqref{eq:prod_expanded} for $u$ between $1$ and $2$, we obtain
\begin{equation}\label{eq:off_diag_terms}
\sideset{}{'}\sum_{\substack{\mathbf{p}_1 \in \mathcal{P}_1 ,
\dots , \mathbf{p}_k \in \mathcal{P}_k \\ \mathbf{q}_1 \in
\mathcal{Q}_1 , \dots , \mathbf{q}_k \in \mathcal{Q}_k}}
\frac{\exp\left(\ii 2 F(\mathbf{p} , \mathbf{q}) \right) -
\exp\left(\ii F(\mathbf{p} , \mathbf{q}) \right)}{\ii F(\mathbf{p}
, \mathbf{q}) \sqrt{\prod_{j=1}^k \left( \prod_{\ell_j=1}^{m_j}
p_{j,\ell_j} \right) \left(\prod_{\ell_j=1}^{n_j}
q_{j,\ell_j}\right)}}
\end{equation}
where $\sum'$ denotes that we are summing only over the
non-diagonal terms, those terms where $F(\mathbf{p} , \mathbf{q})
\neq 0$.

Recall that, without loss of generality, we assumed
$\lambda_1>\dots>\lambda_k$. Assume
\begin{equation}\label{eq:cond_on_non_diag}
\log\left(\frac{q_{1,1}\dots q_{1,n_1}}{p_{1,1}\dots
p_{1,m_1}}\right) \neq 0
\end{equation}
Since  $p_{1,\ell} \leq \exp\left(\frac{N^{\lambda_1}}{40 k
m_1}\right)$ and $q_{1,\ell} \leq \exp\left(\frac{N^{\lambda_1}}{40
k n_1}\right)$, we have
\begin{equation*}
\exp\left(N^{\lambda_1}\right) \left|\log\left(\frac{q_{1,1}\dots
q_{1,n_1}}{p_{1,1}\dots p_{1,m_1}}\right)\right| > \frac12
\exp\left((1-\frac{1}{40k}) N^{\lambda_1}\right)
\end{equation*}
which follows from the fact that if $m,n$ are positive integers,
and $m\neq n$, then $|\log(m/n)| > 1/(2\min(m,n))$. Furthermore,
for $j>1$, for $\mathbf{p}_j \in \mathcal{P}_j$ and $\mathbf{q}_j
\in \mathcal{Q}_j$, then for sufficiently large $N$,
\begin{equation*}
\exp\left(N^{\lambda_j}\right) \left|\log\left(\frac{q_{j,1}\dots
q_{j,n_j}}{p_{j,1}\dots p_{j,m_j}}\right)\right| \leq \frac{1}{40
k} \exp\left(N^{\lambda_j}\right) N^{\lambda_j} <
\exp\left(2N^{\lambda_j}\right)
\end{equation*}
and so we can conclude that if \eqref{eq:cond_on_non_diag} holds,
\begin{align*}
|F(\mathbf{p} , \mathbf{q})| &>
\exp\left(\left(1-\frac{1}{40k}\right) N^{\lambda_1}\right) -
\sum_{j=2}^k
\exp\left(2N^{\lambda_j}\right) \\
&>\exp\left(\left(1-\frac{1}{20k}\right) N^{\lambda_1}\right)
\end{align*}
for sufficiently large $N$, since $\lambda_j < \lambda_1$ for all
$j>1$.

The contribution of such terms to \eqref{eq:off_diag_terms} is
clearly bounded by
\begin{multline*}
\frac{1}{\exp\left((1-\frac{1}{20k}) N^{\lambda_1}\right)}
\sum_{\substack{\mathbf{p}_1 \in
\mathcal{P}_1 , \dots , \mathbf{p}_k \in \mathcal{P}_k \\
\mathbf{q}_1 \in \mathcal{Q}_1 , \dots , \mathbf{q}_k \in
\mathcal{Q}_k}} \frac{1}{\sqrt{\prod_{j=1}^k \left(
\prod_{\ell_j=1}^{m_j} p_{j,\ell_j} \right)
\left(\prod_{\ell_j=1}^{n_j} q_{j,\ell_j}\right)}} \\
\leq \exp\left(-(1-\frac{1}{20k})N^{\lambda_1}\right)
\exp\left(\sum_{j=1}^k \frac{1}{40 k} N^{\lambda_j} \right) \leq
\exp\left(-(1-\frac{1}{10k}) N^{\lambda_1}\right)
\end{multline*}
once more using the fact that
\begin{equation*}
\sum_{j=2}^k \frac{1}{40k} N^{\lambda_j}  \leq
\frac{1}{20k}N^{\lambda_1}
\end{equation*}
for sufficiently large $N$.

Therefore, as $N$ tends to infinity, we see that terms which
satisfy \eqref{eq:cond_on_non_diag} vanish. Thus, for a non-zero
result in the limit, we must have
\begin{equation*}
\log\left(\frac{q_{1,1}\dots q_{1,n_1}}{p_{1,1}\dots
p_{1,m_1}}\right) = 0
\end{equation*}
That is, we must have $(\mathbf{p}_1,\mathbf{q}_1) \in
\mathcal{D}_1$, where $\mathcal{D}_1$ is defined in
\eqref{eq:defn_D}.

The terms which might possibly contribute to
\eqref{eq:off_diag_terms} are
\begin{equation*}
\sideset{}{'}\sum_{\substack{(\mathbf{p_1} , \mathbf{q}_2) \in
\mathcal{D}_1 \\ \mathbf{p}_2 \in
\mathcal{P}_2 , \dots , \mathbf{p}_k \in \mathcal{P}_k \\
\mathbf{q}_2 \in \mathcal{Q}_2 , \dots , \mathbf{q}_k \in
\mathcal{Q}_k}} \frac{\exp\left(\ii 2 F(\mathbf{p} , \mathbf{q})
\right) - \exp\left(\ii F(\mathbf{p} , \mathbf{q}) \right)}{\ii
F(\mathbf{p} , \mathbf{q}) \sqrt{\prod_{j=1}^k \left(
\prod_{\ell_j=1}^{m_j} p_{j,\ell_j} \right)
\left(\prod_{\ell_j=1}^{n_j} q_{j,\ell_j}\right)}}
\end{equation*}
The same argument as above, shows that the terms with
\begin{equation*}
\log\left(\frac{q_{2,1}\dots q_{2,n_2}}{p_{2,1}\dots
p_{2,m_2}}\right) \neq 0
\end{equation*}
contribute
\begin{multline*}
\frac{1}{\exp\left((1-\frac{1}{20k}) N^{\lambda_2}\right)}
\sum_{\substack{(\mathbf{p_1} , \mathbf{q}_2) \in \mathcal{D}_1 \\
\mathbf{p}_2 \in
\mathcal{P}_2 , \dots , \mathbf{p}_k \in \mathcal{P}_k \\
\mathbf{q}_2 \in \mathcal{Q}_2 , \dots , \mathbf{q}_k \in
\mathcal{Q}_k}}  \frac{1}{\sqrt{\prod_{j=1}^k \left(
\prod_{\ell_j=1}^{m_j} p_{j,\ell_j} \right)
\left(\prod_{\ell_j=1}^{n_j} q_{j,\ell_j}\right)}} \\
\leq \exp\left(-(1-\frac{1}{10k}) N^{\lambda_2}\right)
\sum_{(\mathbf{p_1} , \mathbf{q}_2) \in \mathcal{D}_1 }
\frac{1}{\sqrt{\left( \prod_{\ell=1}^{m_1} p_{1,\ell} \right)
\left(\prod_{\ell=1}^{n_1} q_{1,\ell}\right)}}
\end{multline*}
If $(\mathbf{p_1} , \mathbf{q}_2) \in \mathcal{D}_1$, then
\eqref{eq:diag_terms} shows that
\begin{equation*}
\sum_{(\mathbf{p_1} , \mathbf{q}_2) \in \mathcal{D}_1 }
\frac{1}{\sqrt{\left( \prod_{\ell=1}^{m_1} p_{1,\ell} \right)
\left(\prod_{\ell=1}^{n_1} q_{1,\ell}\right)}} \ll (\log N)^{n_j}
\end{equation*}
which is negligible compared to $\exp\left(-(1-\frac{1}{10k})
N^{\lambda_2}\right)$. Hence, finally, these terms do not
contribute.

Repeating the argument for $\lambda_j$, $j=3,4,\dots, k$ we see
that any term with
\begin{equation*}
\log\left(\frac{q_{j,1}\dots q_{j,n_j}}{p_{j,1}\dots
p_{j,m_j}}\right) \neq 0
\end{equation*}
has a vanishing contribution to the large-$N$ limit. Therefore,
the main term must come from those terms for which
\begin{equation*}
\log\left(\frac{q_{j,1}\dots q_{j,n_j}}{p_{j,1}\dots
p_{j,m_j}}\right) = 0
\end{equation*}
for all $j$. Such terms are the diagonal terms, and their
contribution has been calculated above. This completes the proof
of Lemma~\ref{lem:mmts_primes}.
\end{proof}

\begin{proof}[Proof of Theorem~\ref{lem:mmts_zeta}]
Recall that
\begin{equation*}
L_\lambda(N,u) = \frac{\log\zeta(\tfrac12 + \ii u
e^{N^\lambda})}{\sqrt{\log N}} .
\end{equation*}
and
\begin{equation*}
P(\lambda,n;k,N,u) = \frac{1}{\sqrt{\log N}} \sum_{p \leq
\exp\left(\frac{N^{\lambda_\ell}}{40 k m_\ell}\right)}
\frac{p^{-\ii u e^{N^\lambda}}}{\sqrt{p}}
\end{equation*}
Let
\begin{align*}
\epsilon(\lambda,n) &= \epsilon(\lambda,n;k,N,u) \\
&= L_\lambda(N,u) - P(\lambda,n;k,N,u)
\end{align*}
so that, if we write $T = \exp(N^{\lambda_j})$, then changing
variables to $t = T u$,
\begin{multline}\label{eq:apply_Selberg}
\int_1^2 \left|\epsilon(\lambda_j,m_j) \right|^{2 k m_j} \;\dd u \\
= \frac{1}{(\log N)^{k m_j}} \frac{1}{T} \int_T^{2T} \left|
\log\zeta(\tfrac12 + \ii t) - \sum_{p \leq T^{1/40 k m_\ell}}
\frac{p^{-\ii t}}{\sqrt{p}}
\right|^{2 k m_j} \;\dd t \\
= O\left(\frac{(k m_j)^{4 k m_j} e^{A k m_j}}{(\log N)^{k m_j}}
\right)
\end{multline}
by Theorem~\ref{thm:Selberg}. Since the $m_j$ are fixed, this
tends to zero as $N\to\infty$.

Consider
\begin{equation}\label{eq:int_difference}
\int_1^2 \left|\prod_{j=1}^k L_{\lambda_j}(N,u)^{m_j}
\overline{L_{\lambda_j}(N,u)}^{n_j} - \prod_{j=1}^k
P(\lambda_j,m_j)^{m_j} \overline{P(\lambda_j,n_j) }^{n_j} \right|
\;\dd u
\end{equation}
Writing $L_{\lambda_j}(N,u)$ in terms of $P(\lambda_j,m_j)$ and
$\epsilon(\lambda_j,m_j)$, we see that the term inside the modulus
signs equals
\begin{equation*}
\sum_{\substack{0\leq \alpha_1 \leq m_1 , \dots , 0\leq \alpha_k \leq m_k \\
0 \leq \beta_1 \leq n_1 , \dots, 0 \leq \beta_k \leq n_k \\
\sum \alpha_j + \beta_j \geq 1}} \prod_{j=1}^k
\binom{m_j}{\alpha_j} \binom{n_j}{\beta_j}  P(\lambda_j,m_j)^{m_j
- \alpha_j} \epsilon(\lambda_j,m_j)^{\alpha_j}
\overline{P(\lambda_j,n_j)}^{n_j - \beta_j}
\overline{\epsilon(\lambda_j,n_j)}^{\beta_j}
\end{equation*}
The integral of this in \eqref{eq:int_difference} is clearly
bounded by
\begin{multline}\label{eq:holder}
\sum_{\substack{0\leq \alpha_1 \leq m_1 , \dots , 0\leq \alpha_k \leq m_k \\
0 \leq \beta_1 \leq n_1 , \dots, 0 \leq \beta_k \leq n_k \\
\sum \alpha_j + \beta_j \geq 1}} \left\{ \prod_{j=1}^k
\binom{m_j}{\alpha_j} \binom{n_j}{\beta_j} \right\}
\\
\times \int_1^2 \prod_{j=1}^k \left|P(\lambda_j,m_j)\right|^{m_j -
\alpha_j} \left|\epsilon(\lambda_j,m_j)\right|^{\alpha_j}
\left|\overline{P(\lambda_j,n_j)}\right|^{n_j - \beta_j}
\left|\overline{\epsilon(\lambda_j,n_j)}\right|^{\beta_j} \;\dd u
\end{multline}
A version of the generalized H\"older inequality states that
\begin{multline*}
\int_1^2 \prod_{j=1}^k |A_j| |B_j| |C_j| |D_j| \;\dd u  \\
\leq \prod_{j=1}^k \left(\int_1^2 |A_j|^{2 k r_j} \;\dd u
\right)^{1/(2 k r_j)} \left(\int_1^2 |B_j|^{2 k s_j} \;\dd u
\right)^{1/(2k s_j)} \\
\times \left(\int_1^2 |C_j|^{2 k t_j} \;\dd u \right)^{1/(2k t_j)}
\left(\int_1^2 |D_j|^{2 k u_j} \;\dd u \right)^{1/(2k u_j)}
\end{multline*}
so long as $\frac{1}{r_j} + \frac{1}{s_j} = 1$ and $\frac{1}{t_j}
+ \frac{1}{u_j} = 1$ for all $j=1,\dots,k$.

Choosing $r_j = m_j/(m_j - \alpha_j)$ and $s_j = m_j / \alpha_j$,
and $t_j = n_j/(n_j-\beta_j)$ and $u_j = n_j / \beta_j$,  we see
that we may bound the above integral by
\begin{multline*}
\prod_{j=1}^k \left(\int_1^2 \left|P(\lambda_j,m_j)\right|^{2 k
m_j} \;\dd u \right)^{\frac{m_j - \alpha_j}{2km_j}}
\left(\int_1^2 \left|\epsilon(\lambda_j,m_j)\right|^{2 k m_j}
\dd u\right)^{\frac{\alpha_j}{2km_j}} \\
\times \left(\int_1^2 \left|\overline{P(\lambda_j,n_j)}\right|^{2
k n_j} \;\dd u \right)^{\frac{n_j - \beta_j}{2kn_j}}
\left(\int_1^2 \left|\overline{\epsilon(\lambda_j,n_j)}\right|^{2
k n_j} \;\dd u \right)^{\frac{\beta_j}{2kn_j}}
\end{multline*}

From \eqref{eq:apply_Selberg}, if $\alpha_j \neq 0$,
\begin{equation*}
\lim_{N\to\infty} \left(\int_1^2
\left|\epsilon(\lambda_j,m_j)\right|^{2 k m_j} \dd u
\right)^{\frac{\alpha_j}{2km_j}} = 0
\end{equation*}
and from Lemma~\ref{lem:mmts_primes} we have
\begin{equation*}
\left(\int_1^2 \left|P(\lambda_j,m_j)\right|^{2 k m_j} \;\dd u
\right)^{\frac{m_j - \alpha_j}{2km_j}} \ll 1
\end{equation*}
Since the sum in \eqref{eq:holder} is over those $\alpha_j$ and
$\beta_j$ such that $\sum \alpha_j + \beta_j \geq 1$, there must be
at least one $j$ with a non-zero $\alpha_j$ or $\beta_j$. Hence, as
$N\to\infty$, all the terms in \eqref{eq:holder} tend to zero. The
sum is over a finite number of terms, so we may conclude that
\begin{equation*}
\lim_{N\to\infty} \int_1^2 \left|\prod_{j=1}^k
L_{\lambda_j}(N,u)^{m_j} \overline{L_{\lambda_j}(N,u)}^{n_j} -
\prod_{j=1}^k P(\lambda_j,m_j)^{m_j} \overline{P(\lambda_j,n_j)
}^{n_j} \right| \;\dd u = 0
\end{equation*}
which implies
\begin{equation*}
\lim_{N\to\infty} \int_1^2 \prod_{j=1}^k L_{\lambda_j}(N,u)^{m_j}
\overline{L_{\lambda_j}(N,u)}^{n_j} \;\dd u \\
= \lim_{N\to\infty} \prod_{j=1}^k P(\lambda_j,m_j)^{m_j}
\overline{P(\lambda_j,n_j) }^{n_j}
\end{equation*}
assuming the limits make sense. Therefore
Theorem~\ref{lem:mmts_zeta} follows from
Lemma~\ref{lem:mmts_primes}.
\end{proof}

\end{document}